\documentclass[12pt,reqno]{amsart}
\usepackage{amsmath,amssymb,amscd,latexsym,amsthm,mathrsfs,comment}
\usepackage[unicode]{hyperref}
\usepackage{hypbmsec}
\newtheorem{Thm}{Theorem}

\newtheorem{Lem}{Lemma}

\theoremstyle{remark}

\begin{document}

\title{A top hat for Moser's four mathemagical rabbits}
\author{Pieter Moree} 
\address{Max-Planck-Institut f\"ur Mathematik, Vivatsgasse 7, D-53111 Bonn, Germany}
\email{moree@mpim-bonn.mpg.de}

\maketitle
\begin{abstract}
If the equation  $1^k+2^k+\cdots+(m-2)^k+(m-1)^k=m^k$ has a solution with $k\ge 2$,
then $m>10^{10^6}$. Leo Moser showed this in 1953 by remarkably
elementary methods. His proof rests on four identities he derives separately.
It is shown here that Moser's result can be derived from a 
von Staudt-Clausen type theorem (an easy proof of which is also presented here). 
In this approach the four identities can be derived uniformly. The mathematical
arguments used
in the proofs were already available during the lifetime of Lagrange (1736-1813).
\end{abstract}

\section{Introduction}
\noindent Consider the Diophantine equation
\begin{equation}
\label{EME}
1^k+2^k+\cdots+(m-2)^k+(m-1)^k=m^k,
\end{equation}
to be solved in integers $(m,k)$ with $m\ge 2$ and $k\ge 1$.
Note that in case
$k=1$ the left-hand side of (\ref{EME}) equals $m(m-1)/2$, and this leads to the (unique)
solution $1+2=3$.
{}From now on
we will assume that $k\ge 2$.
Conjecturally solutions with $k\ge 2$ do not exist (this conjecture was formulated
around 1950 by Paul Erd\H{o}s in a letter to Leo Moser). 
Leo Moser \cite{Moser}
established the following theorem in 1953. 
\begin{Thm} 
\label{Leo}
{\rm (Leo Moser, 1953)}.
If
$(m,k)$ is a solution of \eqref{EME} with $k\ge 2$, then
$m>10^{10^6}$.
\end{Thm}
His result has since then been improved upon.
Butske et al.~\cite{graphs} have shown, by computing rather
than estimating certain quantities in Moser's original proof, that $m>1.485\cdot 10^{9321155}$.
By proceeding along these lines this bound cannot be improved upon substantially.
 Butske et al.~\cite[p. 411]{graphs} expressed
the hope that new insights will eventually make it possible to reach the more natural benchmark
$10^{10^7}$. This hope was recently fulfilled by Gallot, the author, and Zudilin \cite{GMZ},
who showed that $2k/(2m-3)$ must be a convergent of $\log 2$ and made an extensive continued 
fraction computation of $(\log 2)/2N$, with $N$ an appropriate integer, in order to
establish Theorem \ref{GMZ}. Note that their result goes well beyond establishing the
benchmark.
Their approach 
only works for those $N$ for which it can be shown that $N|k$. In \cite{MRU} it was, e.g., shown
that lcm$(1,2,\ldots,200)|k$.
\begin{Thm}
\label{GMZ}
If
$(m,k)$ is a solution of \eqref{EME} with $k\ge 2$, then
$m>10^{10^9}$.
\end{Thm}
Moser's proof of Theorem \ref{Leo} is quite amazing in the sense that he uses only very
elementary number theory. His proof is even mathemagical in the sense that he pulls four rabbits out
of a hat, namely the equations (\ref{M8}), (\ref{M13}), 
(\ref{M16}), and (\ref{M18}), that a solution
$(m,k)$ has to satisfy. He derives each of these equations separately
in a quite ingenious way. In this note we will see that a reproof of Moser's result
can be given, showing that the following result is the top hat the four rabbits were pulled from.\\
\indent Put $S_r(y)=\sum_{j=1}^{y-1}j^r$.
\begin{Thm} {\rm (Carlitz-von Staudt, 1961)}.
\label{thm:clausen}
Let $r$ and $y$ be positive integers. Then
\begin{equation}
\label{staudt}
S_r(y)=\sum_{j=1}^{y-1}j^r\equiv \begin{cases}
0~({\rm mod~}{y(y-1)\over 2}) &{\rm~if~}r{\rm ~is~odd};\cr
-\sum_{(p-1)|r,~p|y}{y\over p}~({\rm mod~}y) &{\rm ~otherwise}.
\end{cases} 
\end{equation}
\end{Thm}
\noindent The latter sum is over the primes $p$ such that both $p-1$ divides $r$ and $p$ divides $y$. (Here and 
in the sequel
the letter $p$ is used to indicate primes.)
Using this result, an easy proof of which will be given in Section \ref{sec:carlitz},
a less mathemagical reproof of Moser's result can be given. For a polished version of 
Moser's original proof, we refer the reader to the extended version of this note \cite{MPRE}.\\
\indent The prime harmonic sum diverges (as Euler already knew) and so given
$\alpha> 1/2$ there exists a largest prime $p(\alpha)$ such that $\sum_{p\le p(\alpha)}1/p< \alpha$.
Moser needed $p(3.16)$ in his proof, but could only estimate it using certain prime number estimates.
His proof is easily adapted to involve $p(3{1\over 6})$ and this was first exactly computed
by Butske et al. \cite{graphs}, leading to an improvement of Moser's bound, namely
\begin{equation}
\label{blobo}
m>\Big(3\prod_{p\le p(3{1\over 6})}p\Big)^{1\over 4}> 1.485\cdot 10^{9321155}.
\end{equation}
\indent We obtain, using Theorem \ref{thm:clausen} and a computer algebra
package like PARI to compute $p(3{1\over 6})=85861889$ and the prime product
in (\ref{blobo}), the following variant of Moser's result.
\begin{Thm}
\label{main}
Suppose that $(m,k)$ is a solution of {\rm (\ref{EME})} with $k\ge 2$. Then\\
{\rm 1)} $m > 1.485\cdot 10^{9321155}$;\\  
{\rm 2)} $k$ is even, $m\equiv 3~({\rm mod~}8)$, $m\equiv \pm 1~({\rm mod~}3)$;\\
{\rm 3)} $m-1$, $(m+1)/2$, $2m-1$, and $2m+1$ are all squarefree;\\
{\rm 4)} if $p$ divides at least one of the above four integers, then $(p-1)|k$;\\
{\rm 5)} the number $(m^2-1)(4m^2-1)/12$ is squarefree and has at least $4990906$ 
prime factors.
\end{Thm}
The proof we give in this note shows that if
Lagrange (1736-1813) had a present-day computer, he could have
proven Theorem \ref{main}.

In order to improve on Theorem \ref{GMZ}  by Moser's approach one needs to find
additional rabbit(s) in the top hat. The interested reader is wished good luck in
finding these elusive animals !

\section{Proof of Theorem {\rm \ref{main}}}
\noindent {\it Proof of Theorem} \ref{main}. We will apply 
Theorem \ref{thm:clausen} with $r=k$.\hfil\break
\indent In case $k$ is odd, we find, on combining (\ref{staudt}) (putting $y=m$) with (\ref{EME}) 
and using the
coprimality of $m$ and $m-1$, that $m=2$ or $m=3$, but these cases are easily excluded. 
(Since $1^k+2^k<(1+2)^k$ for $k>1$, one sees
that $1+2^k=3^k$ has only the solution $k=1$.)
Therefore
$k$ must be even.\\
\indent Take $y=m-1$. Then using (\ref{EME}), the left-hand side of (\ref{staudt}) simplifies to
\begin{equation}
\label{jar}
S_k(m-1)=1^k+2^k+\cdots+(m-2)^k=m^k-(m-1)^k\equiv  1~({\rm mod~}m-1).
\end{equation}
We get from (\ref{staudt}) and (\ref{jar}) that
\begin{equation}
\label{Pi1}
\sum_{p|(m-1),~(p-1)|k}{m-1\over p}+1\equiv 0~({\rm mod~}m-1).
\end{equation}
Suppose there exists $p|(m-1)$ such that $(p-1)\nmid k$. Then on reducing both sides of (\ref{Pi1}) 
modulo $p$ we get $1\equiv 0~({\rm mod~}p)$. This
contradiction shows that in (\ref{Pi1}) the condition $(p-1)|k$ can be dropped, and thus
we obtain 
\begin{equation}
\label{M7}
\sum_{p|(m-1)}{m-1\over p}+1\equiv 0~({\rm mod~}m-1).
\end{equation}
Suppose there exists a prime $p$ dividing $m-1$ such that $p^2$ also divides $m-1$. Then
on reducing both sides modulo $p$, we get $1\equiv 0~({\rm mod~}p)$.
This contradiction shows that $m-1$ must be squarefree. On dividing (\ref{M7}) by $m-1$ 
we obtain 
\begin{equation}
\label{M8}
\sum_{p|(m-1)}{1\over p}+{1\over m-1}\in \Bbb Z.
\end{equation}
\indent Take $y=m$. Then using (\ref{EME}) and $2|k$ we infer from (\ref{staudt}) that
\begin{equation}
\label{five}
\sum_{(p-1)|k,~p|m}{1\over p}\in \Bbb Z.
\end{equation}
Since a sum of reciprocals of distinct primes can never be a positive integer,
we infer that the sum in (\ref{five}) equals zero and hence conclude that if $(p-1)|k$,
then $p\nmid m$. We conclude for example that $(6,m)=1$. Now on considering (\ref{EME}) with
modulus 4 we see that $m\equiv 3~({\rm mod~}8)$.\\
\indent Take $y=m+1$. Then using (\ref{EME}) and the fact that $k$ is even, the left-hand side 
of (\ref{staudt}) simplifies to 
$$S_k(m+1)=S_k(m)+m^k=2(m+1-1)^k\equiv 2~({\rm mod~}m+1).$$
We obtain
\begin{equation}
\label{lila} 
\sum_{p|(m+1),~(p-1)|k}{m+1\over p}+2\equiv 0~({\rm mod~}m+1),
\end{equation}
but by reasoning as in the
case $y=m-1$, it is seen that $p|(m+1)$ implies $(p-1)|k$ and thus 
\begin{equation}
\label{M13}
\sum_{p|(m+1)}{1\over p}+{2\over m+1}\in \Bbb Z.
\end{equation}
{}From (\ref{lila}) and
$m\equiv 3~({\rm mod~}8)$, we derive that $(m+1)/2$ is squarefree.\\
\indent Take $y=2m-1$. On noting that 
$$S_k(2m-1)=\sum_{j=1}^{m-1}(j^k+(2m-1-j)^k)\equiv 2S_k(m)\equiv 2m^k~({\rm mod~}2m-1),$$
we infer that
\begin{equation}
\label{Pi2}
\sum_{p|(2m-1),~(p-1)|k}{2m-1\over p}+2m^k\equiv 0~({\rm mod~}2m-1).
\end{equation}
Since $m$ and $2m-1$ are coprime we infer that if $p|(2m-1)$, then $(p-1)|k$ and $m^k\equiv 1~({\rm mod~}p)$,
and furthermore that $2m-1$ is squarefree. By the
Chinese remainder theorem it then follows that $2m^k\equiv 2~({\rm mod~}2m-1)$,
and hence from (\ref{Pi2}) we obtain 
\begin{equation}
\label{M16}
\sum_{p|(2m-1)}{1\over p}+{2\over 2m-1}\in \Bbb Z.
\end{equation}
\indent Take $y=2m+1$. On noting that
$$S_k(2m+1)=\sum_{j=1}^m(j^k+(2m+1-j)^k)\equiv 2S_k(m+1)\equiv 4m^k~({\rm mod~}2m+1)$$
and proceeding as in the case $y=2m-1$ we obtain 
\begin{equation}
\label{M18}
\sum_{p|(2m+1)}{1\over p}+{4\over 2m+1}\in \Bbb Z.
\end{equation}
We further see that 
$2m+1$ is squarefree.\\
\indent No prime $p>3$ can divide more than one of the integers
$m-1$, $m+1$, $2m-1$, and $2m+1$. Further, 
since $m\equiv 3~({\rm mod~}8)$ and $3\nmid m$,
2 and 3 divide precisely two of these integers.
We infer that $M=(m-1)(m+1)(2m-1)(2m+1)/12$ is a squarefree integer. 
On adding (\ref{M8}), (\ref{M13}), 
(\ref{M16}), and (\ref{M18}), we deduce that
\begin{equation}
\label{M19}
\sum_{p|M}{1\over p}+{1\over m-1}+{2\over m+1}+{2\over 2m-1}+{4\over 2m+1}\ge 4-{1\over 2}-{1\over
3}=3{1\over 6}.
\end{equation}
One checks that the only solutions of (\ref{M8}) with $m\le 1000$ are $3,7,$ and $43$. These
are easily ruled out by (\ref{M13}). Thus (\ref{M19}) yields (with $\alpha=3.16$)
$\sum_{p|M}{1\over p}>\alpha$.
{}From this it follows that if 
\begin{equation}
\label{M20}
\sum_{p\le x}{1\over p}<\alpha, 
\end{equation}
then
$m^4/3>M>\prod_{p\le x}p$ and hence
\begin{equation}
\label{Pbound}
m > 3^{1/4}e^{\theta(x)/4},
\end{equation}
with $\theta(x)=\sum_{p\le x}\log p$, the Chebyshev $\theta$-function. Since
for example (\ref{M20}) is satisfied with $x=1000$, we find that $m>10^{103}$ and
infer from (\ref{M19}) that we can take $\alpha=3{1\over 6}-10^{-100}$ in (\ref{M20}).
Next one computes (using a computer algebra package) the largest prime $p_k$ such 
that $\sum_{p_j\le p_k}{1\over p_j}<3{1\over 6}$, with
$p_1,p_2,\ldots$ the consecutive primes (note
that $p_k=p(3{1\over 6})$). Here one finds that $k=4990906$ and
$$\sum_{i=1}^{4990906}{1\over p_i}=3.1666666588101728584<3{1\over 6}-10^{-9}.$$ 
By direct computation one finds that $\theta(p_k)=8.58510010694053\cdots \times 10^7$.
Using this we infer from (\ref{Pbound}) the inequality (\ref{blobo}), and hence
part 1 of the theorem is proved.\\
\indent Notice that along our way towards proving part 1, the 
remaining parts of the theorem have also been proved. \qed

\section{Proof of the Carlitz-von Staudt theorem}
\label{sec:carlitz}
Carlitz \cite{Carlitz} gave a proof of Theorem \ref{thm:clausen} using finite
differences and stated that the result is due
to von Staudt. When $r$ is odd, he claims that $S_r(y)/y$ is an integer, which
is not always true (it is true though that $2S_r(y)/y$ is always an integer). The 
author \cite{Canada} gave a reproof using the theory of primitive
roots and Kellner \cite{Kellner2} a reproof (for $r$ even only) using Stirling numbers of 
the second kind.
Here a reproof will be given that is easier than all the above. It uses only the following
result of Lagrange.
\begin{Thm}
\label{Lagrange}
If $f$ is a one-variable polynomial of degree $n$ over $\mathbb Z/p\mathbb Z$, then it
cannot have more than $n$ roots unless it is identically zero.
\end{Thm}
\noindent {\it Proof}. See, e.g., the book of Rose \cite[Theorem 2.2, p. 39]{Rose}. \qed
\begin{Lem}
\label{donz}
Suppose that $(p-1)\nmid r$. Then the equation $x^r\not\equiv 1~({\rm mod~}p)$
has a solution.
\end{Lem}
\noindent {\it Proof}. Let $r_1$ be the smallest positive integer such 
that $r_1\equiv r~({\rm mod~}p-1)$. Then $r_1<p-1$. Suppose
that $x^{r}\equiv 1~({\rm mod~}p)$ for every $x\in \{1,2,\ldots,p-1\}$. Then by Fermat's little theorem 
we also have $x^{r_1}\equiv 1~({\rm mod~}p)$ for every $x\in \{1,2,\ldots,p-1\}$, contradicting
Lagrange's theorem. \qed
\begin{Lem}
\label{Smod}
Let $p$ be a prime.
We have $$S_r(p)\equiv \epsilon_r(p)~({\rm mod~}p),$$
where
$$\epsilon_r(p)=
\begin{cases}
-1 & {\rm if~}(p-1)|r;\\
 0 & {\rm otherwise}.
\end{cases}
$$
\end{Lem}
\noindent {\it Proof}. If $p-1$ divides $r$ the result follows by Fermat's little theorem. If $(p-1)\nmid r$,
assume that $S_r(p)\not\equiv 0~({\rm mod~}p)$. Let $a$ be an integer not divisible by $p$.
Multiplication by $a$ permutes the elements of $\mathbb Z/p\mathbb Z$ and hence
$S_r(p)\equiv a^rS_r(p)~({\rm mod~}p)$, from which we infer that $a^r \equiv 1~({\rm mod~}p)$.
Thus $a^r\equiv 1~({\rm mod~}p)$ for $a=1,2,\ldots,p-1$. Invoking Lemma \ref{donz} gives a contradiction, and hence our assumption that
$S_r(p)\not\equiv 0~({\rm mod~}p)$ must have been false. \qed\\

\noindent The usual proof of this result makes use of the existence of a primitive root modulo $p$, which
provides a solution to $x^r\not\equiv 1~({\rm mod~}p)$ in case $(p-1)\nmid r$. The proof given here 
only makes use of the more elementary theorem of Lagrange, Theorem \ref{Lagrange}.

\begin{Lem}
\label{pripower}
If $p$ is odd or $p=2$ and $r$ is even, we have
$S_r(p^{\lambda+1})\equiv pS_r(p^{\lambda})~({\rm mod~}p^{\lambda+1})$.
\end{Lem}
\noindent {\it Proof}. Every $j$ with $0\le j<p^{\lambda+1}$ can be uniquely written as 
$j=\alpha p^{\lambda}+\beta$ with $0\le \alpha <p$ and $0\le \beta<p^{\lambda}$.
Hence we obtain on invoking the binomial theorem that
$$S_r(p^{\lambda+1})=\sum_{\alpha=0}^{p-1}\sum_{\beta=0}^{p^{\lambda}-1}(\alpha p^{\lambda}+\beta)^r
\equiv
p\sum_{\beta=0}^{p^{\lambda}-1}\beta^r+rp^{\lambda}\sum_{\alpha=0}^{p-1}\alpha\sum_{\beta=0}^{p^{\lambda}-1}\beta^{r-1}~({\rm mod~}p^{2\lambda}) .$$
Since the first sum equals $S_r(p^{\lambda})$ and $2\sum_{\alpha=0}^{p-1}\alpha=p(p-1)\equiv 0~({\rm mod~}p)$,
the result follows. \qed\\

\noindent {\it Proof of Theorem} \ref{thm:clausen}.
First let us consider the case where $r$ is odd. 
We proceed by induction on $y$.
Assume $S_r(m)$  is a multiple of $m(m-1)/2$.
We need to show that $S_r(m+1) = S_r(m)+m^r$ is a multiple of $m(m+1)/2$.

If $m$ is even, we have that $m/2$ divides $S_r(m)$. 
But $$S_r(m+1) = (1^r+m^r)+(2^r+(m-1)^r)+\cdots + ( ({m\over 2})^r+({m\over 2}+1)^r),$$ which is a multiple of 
$m+1$ as each pair above is. Thus, $S_r(m+1)$ is a multiple of $m/2$ as well as of $m+1$. Since these are 
coprime $S_r(m+1)$ is a multiple of $m(m+1)/2$.

If $m$ is odd, then $m|S_r(m)$. 
But $$S_r(m+1) = (1^r+m^r)+(2^r+(m-1)^r)+\cdots+ ({m+1\over 2} )^r,$$ which is a multiple of 
$(m+1)/2$ as each term is. 
Thus $S_r(m+1)$ is a multiple of both $m$ and $(m+1)/2$, which are coprime, and hence
it is a multiple of $m(m+1)/2$.

Next we consider the case where $r$ is even. Suppose that $p^f|y$. Then
\begin{equation}
\label{grigor}
S_r(y)=\sum_{\alpha=0}^{{y\over p^f}-1}\sum_{\beta=0}^{p^f-1}(\alpha p^f+\beta)^r\equiv
{y\over p^f}S_r(p^f)~({\rm mod~}p^f).
\end{equation}
By the Chinese remainder theorem it is enough to show that
$$S_r(y)\equiv {y\over p}\epsilon_r(p) ~({\rm mod~}p^{e_p}),$$ where $y=\prod_p p^{e_p}$ is a
factorization of $y$ into prime powers $p^{e_p}$. By (\ref{grigor}), Lemma \ref{pripower},
and Lemma \ref{Smod}, we then infer that
$$S_r(y)\equiv {y\over p^{e_p}}S_r(p^{e_p})\equiv {y\over p}S_r(p)\equiv {y\over p}\epsilon_r(p)
~({\rm mod~}p^{e_p}),$$
thus concluding the proof. \qed

\section{Concluding remarks}
A further application of the Carlitz-von Staudt theorem is to show that Giuga's conjecture
(1950) and Agoh's conjecture (1990) are equivalent; see Kellner \cite{Kellner2}. Giuga's conjecture
states that if $n\ge 2$, then $S_{n-1}(n)\equiv -1~({\rm mod~}n)$ if and only if $n$ is prime. Agoh's conjecture
states  that if $n\ge 2$, then $nB_{n-1}\equiv -1~({\rm mod~}n)$ if and only if $n$ is prime, where $B_r$ 
denotes the $r$th Bernoulli number.\\
\indent The author has generalized the Carlitz-von Staudt theorem to deal with consecutive $r$th
powers in arithmetic progression; see \cite{Canada}. However, the method of proof in case $r$ is odd given in 
Section \ref{sec:carlitz} no longer applies in this more general situation.\\
\indent That Theorem \ref{thm:clausen} can be used to reprove Moser's result
was first observed by the author in \cite{Oz}, where it played a key role in the study of the more general
equation $1^k+2^k+\cdots+(m-1)^k=am^k$. The presentation given here also draws on
computer improvements since 1996 and \cite{graphs}. The proof of Theorem 
\ref{thm:clausen} given here is clearly easier than those given in \cite{Carlitz, Kellner2, Canada}, and
is the main new contribution in this note.\\
\indent Some variants of the Erd\H{o}s-Moser problem 
require computing $p(\alpha)$ for $\alpha>3{1\over 6}$; see, e.g., \cite{Oz}
The largest value for which $p(\alpha)$ has 
been computed is $\alpha=4$. Bach et al.~\cite{BKS} found that $p(4)=1801241230056600467$, but whereas
the computation of $p(3{1\over 6})$ is straightforward with a computer algebra package, 
computation of $p(4)$ is rather more involved (using the
Meissel-Lehmer algorithm). For $\alpha>4$, one
presently has to resort to deriving a sharp lower bound for $p(\alpha)$ and here one is
forced, as was Moser, to use prime number estimates; cf. \cite{Oz}.\\

\noindent {\tt Acknowledgement}. The argument for Theorem \ref{thm:clausen} in case $r$ is odd
was suggested to me by B.~Sury and the proof of Lemma \ref{Smod} by D.~Zagier. I would like
to thank W.~Moree, J.~Sondow, and the referee for comments on an earlier version.


\begin{thebibliography}{99}
\bibitem{BKS} E. Bach, D. Klyve, and J. P. Sorenson, Computing prime harmonic sums. 
\emph{Math. Comp.}  {\bf 78}  (2009) 2283-2305.
\bibitem{graphs} W. Butske, L. M. Jaje, and D. R. Mayernik, On the 
equation $\sum\sb {p\vert N}\frac1p+\frac1N=1$, pseudoperfect numbers, and 
perfectly weighted graphs. \emph{Math. Comp.} {\bf 69} (2000) 407-420.
\bibitem{Carlitz} L. Carlitz, The Staudt-Clausen theorem. \emph{Math. Mag.} 
{\bf 34} (1960/1961) 131--146.  
\bibitem{GMZ} Y. Gallot, P. Moree, and W. Zudilin, The Erd\H{o}s-Moser equation
$1^k+2^k+\cdots+(m-1)^k=m^k$ revisited using continued fractions (2009), available at 
{\tt http://front.math.ucdavis.edu/0907.1356}. 
\bibitem{Kellner2} B. C. Kellner, The equivalence of Giuga's and Agoh's conjectures (2004), 
available at {\tt http://front.math.ucdavis.edu/0409.5259}.
\bibitem{Canada} P. Moree, On a theorem of Carlitz-von Staudt. 
\emph{C. R. Math. Rep. Acad. Sci. Canada}  {\bf 16}  (1994) 166-170.
\bibitem{Oz} \rule[4pt]{3em}{1.2pt}, Diophantine equations of Erd\H{o}s-Moser type. 
{\it Bull. Austral. Math. Soc.}  {\bf 53} (1996) 281-292.
\bibitem{MPRE} \rule[4pt]{3em}{1.2pt}, Moser's mathemagical work on the 
equation $1^k+2^k+\dots+(m-1)^k=m^k$ (to appear), available at
{\tt http://www.mpim-bonn.mpg.de/preprints/retrieve}.
\bibitem{MRU}  P. Moree, H. J. J.  te Riele, and J. Urbanowicz, Divisibility 
properties of integers $x, k$ satisfying 
$1\sp k+\cdots+(x-1)\sp k=x\sp k$. \emph{Math. Comp.}  {\bf 63}  (1994) 799-815.
\bibitem{Moser} L. Moser, On the 
diophantine equation $1\sp n+2\sp n+3\sp n+\cdots +(m-1)\sp n=m\sp n$.
\emph{Scripta Math.} {\bf 19} (1953) 84-88.
\bibitem{Rose} H. E. Rose, \emph{A Course in Number Theory,} 
Oxford University Press, New York, 1988.

\end{thebibliography}
\end{document}